\documentclass[11pt]{article}
\usepackage{epic,latexsym,amssymb}
\usepackage{color}
\usepackage{tikz}
\usepackage{graphicx}
\usepackage{geometry,verbatim,amsmath}
\usetikzlibrary{calc}

\newtheorem{thm}{Theorem}[section]

\newtheorem{claim}[thm]{Claim}
\newtheorem{cor}[thm]{Corollary}
\newtheorem{obser}[thm]{Observation}
\newtheorem{prop}[thm]{Proposition}
\newtheorem{conj}[thm]{Conjecture}
\newtheorem{problem}[thm]{Problem}

\newcommand{\proof}{\noindent{\bf Proof.\ }}
\newcommand{\qed}{\hfill $\square$ \medskip}

\newcommand{\smallqed}{{\tiny ($\Box$)}}
\newcommand{\gt}{\gamma_t}
\newcommand{\cartprod}{\,\square\,}
\newcommand{\barx}{\overline{x}}
\newcommand{\barv}{\overline{v}}

\newcommand{\barw}{\overline{w}}

\newcommand{\gpr}{\gamma_{\rm pr}}
\newcommand{\tr}{\gamma_{\rm tr}}

\newenvironment{unnumbered}[1]{\trivlist \item [\hskip \labelsep {\bf
#1}]\ignorespaces\it}{\endtrivlist}

\let\oldenumerate\enumerate
\renewcommand{\enumerate}{
  \oldenumerate
  \setlength{\itemsep}{0pt}
  \setlength{\parskip}{0pt}
  \setlength{\parsep}{0pt}
}

\def \HG {H_{_G}}
\def \HGc {H_{_G}^{c}}

\def \HX {H_{_X}}
\def \HY {H_{_Y}}

\def \eX {e_{_X}}
\def \eY {e_{_Y}}

\def \TX {T_{_X}}
\def \TXc {T_{_X}^c}

\begin{document}

\title{(Total) Domination in Prisms}

\author{
Jernej Azarija $^{a,b}$
\and
Michael A. Henning $^{c}$
\and
Sandi Klav\v zar $^{a,b,d}$
}

\date{}

\maketitle

\begin{center}
$^a$ Faculty of Mathematics and Physics, University of Ljubljana, Slovenia \\
jernej.azarija@gmail.com \\
sandi.klavzar@fmf.uni-lj.si \\
\medskip

$^b$ Institute of Mathematics, Physics and Mechanics, Ljubljana, Slovenia\\
\medskip

$^c$ Department of Pure and Applied Mathematics \\
University of Johannesburg, South Africa \\
mahenning@uj.ac.za \\
\medskip

$^d$ Faculty of Natural Sciences and Mathematics, University of Maribor, Slovenia\\
\end{center}

\date{}
\maketitle

\begin{abstract}
With the aid of hypergraph transversals it is proved that $\gt(Q_{n+1}) = 2\gamma(Q_n)$, where $\gt(G)$ and $\gamma(G)$ denote the total domination number and the domination number of $G$, respectively, and $Q_n$ is the $n$-dimensional hypercube. More generally, it is shown that if $G$ is a bipartite graph, then $\gt(G \cartprod K_2) = 2\gamma(G)$. Further, we show that the bipartite condition is essential by constructing, for any $k \ge 1$, a (non-bipartite) graph $G$ such that $\gamma_t(G\cartprod K_2) = 2\gamma(G) - k$. Along the way several domination-type identities for hypercubes are also obtained.
\end{abstract}


\medskip
\noindent{\bf AMS Subj. Class.:} 05C69, 91C65, 05C76, 94B65

\section{Introduction}

Domination and total domination in graphs are very well studied in the literature, here we study these concepts in prisms of graphs, in particular in hypercubes. To determine the domination number $\gamma$ of the $n$-dimensional hypercube $Q_n$, is a fundamental problem in coding theory, 
computer science, and of course in graph theory. In coding theory,
the problem equivalent to the determination of $\gamma(Q_n)$ is to
find the size of a minimal covering code of length $n$ and covering
radius $1$. In computer science, different distribution type problems
on interconnection networks can be modelled by domination invariants,
where hypercubes in turn form a central model for interconnection
networks.

To determine $\gamma(Q_n)$ turns out to be an intrinsically difficult problem. To date, exact values are only known for $n \le 9$. These results are summarized in Table~\ref{table:hypercubes}.

\begin{table}[ht!]
\label{table:hypercubes}
\begin{center}
\begin{tabular}{|c|cccccccccc|}
	\hline
	\textit{n} & 1 & 2 & 3 & 4  & 5  & 6  & 7  & 8  & 9 & 10 \\
	\hline\hline
	$\gamma(Q_n)$ & 1 & 2 & 2 & 4 & 7 & 12 & 16 & 32 & 62 & 107-120 \\
	\hline
\end{tabular}
\end{center}
\caption{Domination numbers of hypercubes up to dimension $10$}
\end{table}

We have checked these values by formulating an integer linear program and solving it with CPLEX. The result $\gamma(Q_9) = 62$ due to {\"O}sterg{\aa}rd and Blass~\cite{ostergard-2001} actually presented a breakthrough back in 2001. The value of $\gamma(Q_{10})$ is currently unknown, see~\cite{bertolo-2004} for the present best lower bound as given in Table~\ref{table:hypercubes} and~\cite{keri-2005} for the present best upper bound. 

Total domination $\gamma_t$ is, besides classical domination, among the most fundamental concepts in domination theory. It has in particular been extensively investigated on Cartesian product graphs (cf.\ \cite{choudhary-2013, henning-2005, lu-2010}), which was in a great part motivated by the famous Vizing's conjecture~\cite{bresar-2012}. Specifically, $\gamma_t(Q_n)$ was recently investigated in the thesis~\cite{verstraten-2014} under the notion of a {\em binary covering code of empty spheres of length $n$ and radius $1$}. In particular, values $\gamma_t(Q_n)$ for $n \leq 10$ were computed and some bounds established. These exact values intrigued us to wonder whether there exists some general relation between the domination number and the total domination number in hypercubes.  

From our perspective it is utmost important that $Q_n$ can be represented as the $n^{\rm th}$ power of $K_2$ with respect to the Cartesian product operation $\cartprod$, that is, $Q_1 =K_2$ and $Q_n = Q_{n-1} \cartprod K_2$ for $n \ge 2$. Our immediate aim in this paper is to prove that $\gt(Q_{n+1}) = 2\gamma(Q_n)$ holds for all $n \ge 1$. For this purpose, we prove the following much more general result that the total domination of a bipartite prism of a graph $G$ is equal to twice the domination number of $G$.

\begin{thm}
\label{t:prism}
If $G$ is a bipartite graph, then
$$\gt(G \cartprod K_2) = 2\gamma(G)\,.$$
\end{thm}

Since $Q_n$, $n \ge 1$, is a bipartite graph, as a special case of Theorem~\ref{t:prism} we note that $\gt(Q_{n+1}) = 2\gamma(Q_n)$. Our second aim is to show that the bipartite condition in the statement of Theorem~\ref{t:prism} is essential. For this purpose, we prove the following result.

\begin{thm}
\label{t:nonbipartite}
For each integer $k \ge 1$, there exists a connected graph $G_k$ satisfying
\[
\gt(G_k \cartprod K_2) - 2\gamma(G_k) = k.
\]
\end{thm}

We proceed as follows. In the next section concepts used throughout the paper are introduced and known facts and results needed are recalled. In particular, the state of the art on $\gamma(Q_n)$ is surveyed. In Section~\ref{sec:proof-of-main}, Theorem~\ref{t:prism} is proved and several of its consequences listed. A proof of Theorem~\ref{t:nonbipartite} is given in Section~\ref{sec:proof-of-non-bipartite}. We conclude the paper with some open problems. In particular we conjecture that the equality in Theorem~\ref{t:prism} holds for almost all graphs.

\section{Preliminaries}

Let $G$ be a graph with vertex set $V(G)$ and edge set $E(G)$. The \emph{order} of $G$ is denoted by $n(G) = |V(G)|$. The \emph{open neighborhood} of a vertex $v$ in $G$ is $N_G(v) = \{u \in V(G) \, | \, uv \in E(G)\}$ and the \emph{closed neighborhood of $v$} is $N_G[v] = \{v\} \cup N_G(v)$. 

For graphs $G$ and $H$, the \emph{Cartesian product} $G \cartprod H$ is the graph with vertex set $V(G) \times V(H)$ where vertices $(u_1,v_1)$ and $(u_2,v_2)$ are adjacent if and only if either $u_1 = u_2$ and $v_1v_2 \in E(H)$ or $v_1 = v_2$ and $u_1u_2 \in E(G)$. If $(u,v)\in V(G\cartprod H)$, then the subgraph of $G\cartprod H$ induced by the vertices of the form $(u,x)$,  $x\in V(H)$, is isomorphic to $H$; it is called the {\em $H$-layer} (through $(u,v)$). Analogously $G$-layers are defined.
The \emph{prism} of a graph $G$ is the graph $G \cartprod K_2$. Note that $G \cartprod K_2$ contains precisely two $G$-layers. Further, if $G$ is a bipartite graph, then we call the prism $G \cartprod K_2$ the \emph{bipartite prism} of $G$. As already mentioned in the introduction, $Q_n$ is a (bipartite) prism because $Q_n = Q_{n-1}\cartprod K_2$.

A \emph{dominating set} of a graph $G$ is a set $S$ of vertices of $G$ such that every vertex in $V(G) \setminus S$ is adjacent to at least one vertex in $S$, while a \emph{total dominating set} of $G$ is a set $S$ of vertices of $G$ such that every vertex in $V(G)$ is adjacent to at least one vertex in $S$.  The \emph{domination number} of $G$, denoted by $\gamma(G)$, is the minimum cardinality of a dominating set of $G$ and the \emph{total domination number} of $G$, denoted by $\gt(G)$, is the minimum cardinality of a total dominating set of $G$. We refer to the books~\cite{haynes-1988, HeYe_book} for more information on the domination number and the total domination number, respectively.

The values $\gamma(Q_7) = 16$ and $\gamma(Q_8) = 32$ also follow from the following result which gives exact values for two infinite families of hypercubes.

\begin{thm}
\label{thm:infinite-families}
If $k\ge 1$, then $\gamma(Q_{2^k-1}) = 2^{2^k-k-1}$ and $\gamma(Q_{2^k}) = 2^{2^k-k}$.
\end{thm}

The first assertion of Theorem~\ref{thm:infinite-families} is based on the fact that hypercubes $Q_{2^k-1}$ contain perfect codes, cf.~\cite{harary-1993}. Since the domination number of a graph with a perfect code is equal to the size of such a code, the assertion follows. Knowing the existence of such codes, by the divisibility condition one immediately infers that $Q_n$ contains a perfect code if and only if $n=2^k-1$ for some $k\ge 1$. Lee~\cite[Theorem 3]{lee-2001} further proved that this is equivalent to the fact that $Q_n$ is a regular covering of the complete graph $K_{n+1}$. The second assertion of Theorem~\ref{thm:infinite-families} is due to van Wee~\cite{wee-1988}. Related aspects of domination in hypercubes were investigated in~\cite{weichsel-1994}.

A set $S$ of vertices in $G$ is a \emph{paired}-\emph{dominating set} if every vertex of $G$ is adjacent to a vertex in $S$ and the subgraph induced by $S$ contains a perfect matching (not necessarily as an induced subgraph). The minimum cardinality of a paired-dominating set of $G$ is the \emph{paired-domination number} of $G$, denoted $\gpr(G)$. A survey on paired-domination in graphs can be found in~\cite{Desor-2014}. By definition every paired-dominating set is a total dominating set, and every total dominating set is a dominating set. Hence we have the following result first observed by Haynes and Slater~\cite{HaSl98}.

\begin{obser}{\rm(\cite{HaSl98})}
\label{ob:2}
For every isolate-free graph $G$, $\gamma(G) \le \gamma_t(G) \le \gpr(G)$.
\end{obser}

A \emph{total restrained dominating set} of $G$ is a total dominating set $S$ of $G$ with the additional property that every vertex outside $S$ has a neighbor outside $S$; that is, $G[V(G) \setminus S]$ contains no isolated vertex. The \emph{total restrained domination number} of $G$, denoted $\tr(G)$, is the minimum cardinality of a total restrained dominating set. The concept of total restrained domination in graphs was introduced by Telle and Proskurowksi~\cite{TePr97} as a vertex partitioning problem.  By definition every total restrained dominating set if a total dominating set, implying the following observation.

\begin{obser}{\rm(\cite{HaSl98})}
\label{ob:3}
For every isolate-free graph $G$, $\gt(G) \le \tr(G)$. 
\end{obser}

The \emph{open neighborhood hypergraph}, abbreviated ONH, of $G$ is the hypergraph $H_G$ with vertex set $V(H_G) = V(G)$ and with edge set $E(H_G) = \{ N_G(x) \mid x \in V(G) \}$ consisting of the open neighborhoods of vertices in $G$. The \emph{closed neighborhood hypergraph}, abbreviated CNH, of $G$ is the hypergraph $H_G^c$ with vertex set $V(H_G^c) = V(G)$ and with edge set $E(H_G) = \{ N_G[x] \mid x \in V(G) \}$ consisting of the closed neighborhoods of vertices in $G$.

A subset $T$ of vertices in a hypergraph $H$ is a \emph{transversal} (also called \emph{vertex cover} or \emph{hitting set}) if $T$ has a nonempty intersection with every edge of $H$. The \emph{transversal number} $\tau(H)$ of $H$ is the minimum size of a transversal in $H$. A transversal of size~$\tau(H)$ is called a $\tau(H)$-set.

The transversal number of the ONH of a graph is precisely the total domination number of the graph, while the transversal number of the CNH of a graph is precisely the domination number of the graph. We state this formally as follows.

\begin{obser}
\label{ob:1}
If $G$ is a graph, then $\gt(G) = \tau(\HG)$ and $\gamma(G) = \tau(\HGc)$.
\end{obser}

We shall also need the following result from~\cite{HeYe08} (see also~\cite{HeYe_book}).

\begin{thm}{\rm (\cite{HeYe08})}
\label{thm_ONH}
The ONH of a connected bipartite graph consists of two components {\rm (}which are induced by the two partite sets of the graph{\rm )}, while the ONH of a connected graph that is not bipartite is connected.
\end{thm}

\section{Proof of Theorem~\ref{t:prism} and its Consequences}
\label{sec:proof-of-main}

In this section, we first present a proof of Theorem~\ref{t:prism}. Recall its statement.

\bigskip
\noindent \textbf{Theorem~\ref{t:prism}} \emph{If $G$ is a bipartite graph, then
$\gt(G \cartprod K_2) = 2\gamma(G)$.}

\bigskip
\proof
Note first that $K_1\cartprod K_2 = K_2$, hence the assertion of the theorem holds for $G=K_1$. Since we can apply the result to each component of the bipartite graph $G$, we may assume that $G$ is connected. Hence in the rest of the proof let $G$ be a connected bipartite graph of order at least $2$.

Let $G_1$ and $G_2$ be the $G$-layers of $G \cartprod K_2$, and let $V_i = V(G_i)$ for $i \in [2]$. For notational convenience, for each vertex $v$ in $G_1$ we denote the corresponding vertex in $G_2$ that is adjacent to $v$ in $G \cartprod K_2$ by $v'$. Thus, the set $\cup_{v \in V_1} \{vv'\}$ of edges between $V_1$ and $V_2$ in $G \cartprod K_2$ forms a perfect matching in $G \cartprod K_2$.

Since $G$ is a bipartite graph, $G\cartprod K_2$ is bipartite as well. Let $X$ and $Y$ be the partite sets of $G \cartprod K_2$. If $w \in \{v,v'\}$ for some vertex $v \in V_1$, then we define the \emph{complement} of the vertex $w$ to be the vertex $\barw \in \{v,v'\} \setminus \{w\}$. We note that if $w \in V_{3-i}$, then $\barw \in V_i$ for $i \in [2]$. Further, we note that $w$ and $\barw$ belong to different partite sets of $G \cartprod K_2$.

Let $H$ be the ONH of $G \cartprod K_2$. By Theorem~\ref{thm_ONH}, $H$ consists of two components that are induced by the two partite sets, $X$ and $Y$, of $G$. Let $\HX$ and $\HY$ be the two components of $H$, where $V(\HX) = X$ and $V(\HY) = Y$. We note that each edge in $\HX$ and $\HY$ corresponds to the open neighborhood of some vertex in $Y$ and some vertex in $X$, respectively, in $G$. For each vertex $w$ in $G \cartprod K_2$, let $e_w$ be the associated hyperedge in $H$; that is, $e_w = N_G(w)$.

We proceed further with the following series of claims. 

\begin{claim}
\label{claim:1}
The hypergraphs $\HX$ and $\HY$ are isomorphic.
\end{claim}

\proof Let $f \colon X \to Y$ be the function that assigns to each vertex $x \in X$ the vertex $\barx \in Y$. Then, $f$ is a bijection between the vertex set of $\HX$ and $\HY$. Suppose that $\eX$ is an edge of $\HX$. Thus, $\eX = e_w$ for some vertex $w \in Y$. The function $f$ maps the edge $\eX$ to the edge $\eY$. We show that $\eY$ is precisely the edge in $\HY$ associated with the vertex $\barw \in X$.

Suppose first that $w \in V_1$. In this case, $\barw = w'$. Let $w$ have degree~$k+1$ in $G \cartprod K_2$, for some $k \ge 1$. Thus, $w$ is adjacent in $G \cartprod K_2$ to $k$ vertices in $V_1$, say to $w_1, w_2, \ldots, w_k$, and to one vertex in $V_2$, namely the vertex $w'$. Since $w \in Y$ and $G \cartprod K_2$ is bipartite, we note that $\{w_1,w_2,\ldots,w_k\} \subseteq V_1 \cap X$ and that $w' \in V_2 \cap X$. Further, the edge $\eX = e_w = \{w_1,w_2,\ldots,w_k,w'\} \in E(\HX)$. Since $f(w_i) = w_i'$ for $i \in [k]$ and $f(w') = w$, the function $f$ maps the edge $\eX$ to the edge $\eY = \{w_1',w_2',\ldots,w_k',w\}$. We note that $\{w_1',w_2',\ldots,w_k',w\} \subseteq Y$, and that $\eY$ is precisely the edge in $\HY$ associated with the vertex $\barw \in X$.

Suppose next that $w \in V_2$. In this case, $w = v'$ for some vertex $v \in V_1$. Thus, $\barw = v$. Let $w$ have degree~$k+1$ in $G \cartprod K_2$, for some $k \ge 1$. Thus, the vertex $v'$ is adjacent in $G \cartprod K_2$ to $k$ vertices in $V_2$, say to $v_1', v_2', \ldots, v_k'$, and to one vertex in $V_1$, namely the vertex $v$. Since $v' \in Y$ and $G \cartprod K_2$ is bipartite, we note that $\{v_1',v_2',\ldots,v_k'\} \subseteq V_2 \cap X$ and that $v \in V_1 \cap X$. Further, the edge $\eX = e_w = \{v_1',v_2',\ldots,v_k',v\} \in E(\HX)$. The function $f$ maps the edge $\eX$ to the edge $\eY = \{v_1,v_2,\ldots,v_k,v'\}$. We note that $\{v_1,v_2,\ldots,v_k,v'\} \subseteq Y$, and that $\eY$ is precisely the edge in $\HY$ associated with the vertex $\barw \in X$.

Suppose that $\eY$ is an edge of $\HY$ and $\eY = e_w$ for some vertex $w \in X$. If the function $f$ maps the edge $\eY$ to $\eX$, then analogously as before, $\eX$ is precisely the edge in $\HX$ associated with the vertex $\barw \in Y$. Thus, the bijective function $f$ preserves adjacency, implying that $\HX$ and $\HY$ are isomorphic.~\smallqed

\begin{claim}
\label{claim:2}
$\gt(G \cartprod K_2) = 2\tau(\HX)$.
\end{claim}

\proof  By Observation~\ref{ob:1}, $\gt(G \cartprod K_2) = \tau(H) = \tau(\HX) + \tau(\HY)$. By Claim~\ref{claim:1}, $\tau(\HX) = \tau(\HY)$, and so $\gt(G \cartprod K_2) = 2\tau(\HX)$.~\smallqed

\begin{claim}
\label{claim:3}
$\gt(G \cartprod K_2) \le 2\gamma(G)$.
\end{claim}

\proof  Let $D$ be a minimum dominating set in $G$, and let $D_1$ and $D_2$ be the copies of $G$ in $G$-layers $G_1$ and $G_2$, respectively. Clearly, $v \in D_1$ if and only if $v' \in D_2$.  The set $D_1 \cup D_2$ is a total dominating set of $G \cartprod K_2$, and so $\gt(G \cartprod K_2) \le |D_1 \cup D_2| = 2|D| = 2\gamma(G)$.~\smallqed

\begin{claim}
\label{claim:4}
$\gamma(G) \le \tau(\HX)$.
\end{claim}
\proof  Let $H^c$ be the CNH of $G$. By Observation~\ref{ob:1}, $\gamma(G) = \tau(H^c)$. We show that $\tau(H^c) \le \tau(\HX)$. Let $\TX$ be a minimum transversal in $\HX$, and so $|\TX| = \tau(\HX)$. We now define the set $\TXc$ as follows. For each vertex $v \in \TX$, we add $v$ to $\TXc$ if $v \in V_1$, otherwise if add $\barv$ to $\TXc$ if $v \in V_2$. We show that $\TXc$ is a transversal in $H^c$. Let $e$ be an arbitrary edge in $H^c$. Thus, $e = N_G[w]$ for some vertex $w$ in $G$. We may assume that the vertices of $G_1$ are named as in the graph $G$, and so $G_1 = G$. In particular, $w \in V_1$. Thus, $\barw = w' \in V_2$.

Suppose that $w \in Y$. In this case, the edge $e_w = N_G(w) = (e \setminus \{w\}) \cup \{\barw\}$ is an edge of $\HX$ and is therefore covered by some vertex, say $z$, of $\TX$. If $z = \barw$, then noting that $\barw \in V_2$, the vertex $w \in \TXc$, and the edge $e$ is therefore covered by a vertex in $\TXc$, namely the vertex $w$. If $z \ne \barw$, then $z$ is a vertex in $e_w$ different from $\barw$. However, $e_w \setminus \{\barw\} = e \setminus \{w\} \subset V_1$, implying that the vertex $z \in V_1$ and therefore $z \in \TXc$. The edge $e$ is therefore covered by a vertex in $\TXc$, namely the vertex $z$. Thus, if $w \in Y$, then the edge $e$ is covered by a vertex in $\TXc$.

Suppose that $w \in X$. We now consider the vertex $\barw \in V_2$. We note that $\barw \in Y$ and that the edge $e_{\barw} = N_G(\barw)$ is an edge of $\HX$. Further, the edge $e_{\barw}$ contains the vertex $w \in V_1$ and all other vertices in $e_{\barw}$ belong to the set $V_2$. Further, if $u$ is a vertex in the edge $e$, then either $u = w$, in which case $u$ also belongs to the edge $e_{\barw}$, or $u \ne w$, in which case $u'$ belongs to the edge $e_{\barw}$. Since the edge $e_{\barw}$ is an edge of $\HX$, it is covered by some vertex, say $z$, of $\TX$. If $z = w$, then noting that $w \in V_1$, the vertex $w \in \TXc$, and the edge $e$ is therefore covered by a vertex in $\TXc$, namely the vertex $w$. If $z \ne \barw$, then $z$ is a vertex in $e_{\barw}$ different from $w$. Thus, $z = u'$ for some vertex $u \in V_1$. Since $u' \in V_2$, the vertex $u \in \TXc$. As observed earlier, $u$ belongs to the edge $e$, implying that the edge $e$ is covered by a vertex in $\TXc$, namely the vertex $u$. Thus, if $w \in X$, then the edge $e$ is covered by a vertex in $\TXc$.

Thus, whenever $w \in X$ or $w \in Y$, the edge $e$ is covered by a vertex in $\TXc$. Since $e$ is an arbitrary edge of $H^c$, this implies that $\TXc$ is a transversal of $H^c$, and therefore that $\tau(H^c) \le |\TXc| = |\TX| = \tau(\HX)$.~\smallqed

\medskip
We now return to the proof of Theorem~\ref{t:prism} one final time. By Claims~\ref{claim:2},~\ref{claim:3}, and~\ref{claim:4}, the following holds.
$$2\tau(\HX) \stackrel{{\rm Claim~\ref{claim:2}}}{=} \gt(G \cartprod K_2) \stackrel{{\rm Claim~\ref{claim:3}}}{\le} 2\gamma(G) \stackrel{{\rm Claim~\ref{claim:4}}}{\le}  2\tau(\HX)\,.$$

Consequently, we must have equality throughout the above inequality chain. In particular, $\gt(G \cartprod K_2) = 2\gamma(G)$. This completes the proof of Theorem~\ref{t:prism}.\qed

\medskip
As an immediate consequence of Theorem~\ref{t:prism} we state that the problems of determining the domination number and the total domination number of hypercubes are equivalent in the following sense:

\begin{cor}
\label{cor:hypercubes}
If $n\ge 1$, then $\gamma_t(Q_{n+1}) = 2\gamma(Q_n)$.
\end{cor}

Combining Corollary~\ref{cor:hypercubes} with Theorem~\ref{thm:infinite-families} we also deduce the following result:

\begin{cor}
\label{cor:infinite-total}
If $k\ge 1$, then $\gamma_t(Q_{2^k+1}) = 2^{2^k-k+1}$ and $\gamma_t(Q_{2^k}) = 2^{2^k-k}$.
\end{cor}

While the first assertion of Corollary~\ref{cor:infinite-total} appears to be new, the second assertion goes back to Johnson~\cite{johnson-1972}, see also~\cite[Theorem 1(b)]{weakley-2006}.

As another consequence of Theorem~\ref{t:prism}, we have the following result.

\begin{cor}
\label{cor:1}
If $G$ is a bipartite graph, then
$$\gt(G \cartprod K_2) = \gpr(G \cartprod K_2) = \tr(G \cartprod K_2)\,.$$
\end{cor}

\proof As shown in the proof of Claim~\ref{claim:3} in Theorem~\ref{t:prism}, if $D_1$ is a minimum dominating set in $G_1$, and $D_2 = \{v' \mid v \in D_1\}$, then the set $D^* = D_1 \cup D_2$ is a total dominating set of $G \cartprod K_2$. We note that $D^*$ is also a paired-dominating set of $G \cartprod K_2$. Further, $|D^*|  = 2\gamma(G)$. By Observation~\ref{ob:2} and Theorem~\ref{t:prism}, this implies that
$$\gt(G \cartprod K_2) \le \gpr(G \cartprod K_2) \le |D^*|  = 2\gamma(G) = \gt(G \cartprod K_2)\,.$$
Consequently, we must have equality throughout the above inequality chain. In particular, $\gt(G \cartprod K_2) = \gpr(G \cartprod K_2)$.
We note that $D^*$ is also a total restrained dominating set of $G \cartprod K_2$. Thus, by Observation~\ref{ob:3}, $\gt(G \cartprod K_2) \le \tr(G \cartprod K_2) \le |D^*|  = 2\gamma(G) = \gt(G \cartprod K_2)$, implying that $\gt(G \cartprod K_2) = \tr(G \cartprod K_2)$.
\qed

As a special case of Theorem~\ref{t:prism} and Corollary~\ref{cor:1}, we have the following result.

\begin{cor}
\label{cor:2}
If $n \ge 1$,  then $\gt(Q_n) = \gpr(Q_{n}) = \tr(Q_{n})$.
\end{cor}

\section{Proof of Theorem~\ref{t:nonbipartite}}
\label{sec:proof-of-non-bipartite}

In this section, we consider general prisms and show that the bipartite condition in the statement of Theorem~\ref{t:prism} is essential. First we recall the trivial lower bound on the total domination number of a graph in terms of the maximum degree of the graph: If $G$ is a graph of order $n$ and maximum degree~$\Delta$ with no isolated vertex, then $\gt(G) \ge n/\Delta$, cf.~\cite[Theorem 2.11]{HeYe_book}.

\begin{prop}
\label{prop:1}
If $k \ge 1$, then $\gt(C_{6k+1}\cartprod K_2) = 2\gamma(C_{6k+1}) - 1$.
\end{prop}

\proof Let $G \cong C_{6k+1}$ for some integer $k \ge 1$. Then, $\gamma(G) = \lceil n(G)/3 \rceil = 2k + 1$. We show that $\gt(G \cartprod K_2) = 4k+1$. Let $G_1$ and $G_2$ be the $G$-layers of $G \cartprod K_2$, where $G_1$ is the cycle $u_1u_2 \ldots u_{6k+1}u_1$ and $G_2$ is the cycle $v_1v_2 \ldots v_{6k+1}v_1$, and where $u_iv_i \in E(G)$. The set
\[
S = \left( \bigcup_{i=0}^{k-1} \{u_{6i+1},u_{6i+2},v_{6k+4},v_{6k+5}\} \right) \cup \{u_{6k+1}\}
\]

\noindent
is a total dominating set of $G \cartprod K_2$, implying that $\gt(G \cartprod K_2) \le |S| = 4k + 1$. Conversely, since $G \cartprod K_2$ is a cubic graph of order~$12k + 2$, the trivial lower bound on the total domination number of $G \cartprod K_2$ is given by $\gt(G \cartprod K_2) \ge (12k+2)/3$, implying that $\gt(G \cartprod K_2) \ge 4k+1$. Consequently, $\gt(G \cartprod K_2) = 4k+1$. As observed earlier, $\gamma(G) = 2k + 1$. Therefore, $\gt(G \cartprod K_2) = 2\gamma(G) - 1$.~\qed

\medskip
We show next that there are connected, non-bipartite graphs $G$ for which the difference $\gt(G \cartprod K_2) - 2\gamma(G)$ can be arbitrarily large. Recall the statement of Theorem~\ref{t:nonbipartite}.

\medskip
\noindent \textbf{Theorem~\ref{t:nonbipartite}} \emph{For each integer $k \ge 1$, there exists a connected graph $G_k$ satisfying
$$\gt(G_k \cartprod K_2) - 2\gamma(G_k) = k.$$
}

\noindent \textbf{Proof.}
For $k = 1$, let $G_1 \cong C_7$. By Proposition~\ref{prop:1}, $\gt(G_1 \cartprod K_2) = 2\gamma(G_1) - 1$. Hence, we assume in what follows that $k \ge 2$. For $i \in [k]$, let $F_i$ be the $5$-cycle $v_{5(i-1)+1} v_{5(i-1)+2} v_{5(i-1)+4} v_{5(i-1)+5} v_{5(i-1)+3} v_{5(i-1)+1}$. Let $G_k$ be obtained from the disjoint union of the cycles $F_1, \ldots, F_k$ by adding the edges $v_{5j}v_{5j+1}$ for $j \in [k-1]$. By construction, $G_k$ is a connected graph of order~$k$. The following two claims determine the domination number of $G_k$ and total domination numbers of the prism $G_k \cartprod K_2$.

\begin{unnumbered}{Claim~A}
For $k \ge 2$, $\gamma(G_k) = 2k$.
\end{unnumbered}
\proof Every dominating set of $G_k$ contains at least two vertices from $V(F_i)$ in order to dominate the vertices in $V(F_i)$ for each $i \in [k]$, and so $\gamma(G_k) \ge 2k$. Conversely, every set consisting of two non-adjacent vertices from each set $V(F_i)$ forms a dominating set of $G_k$, and so $\gamma(G_k) \le 2k$. Consequently, $\gamma(G_k) = 2k$.~\smallqed

\begin{unnumbered}{Claim~B}
For $k \ge 2$, $\gamma_t(G_k \cartprod K_2) = 3k$.
\end{unnumbered}
\proof Let $G_k^1$ and $G_k^2$ be the two copies of the graph $G_k$ in the prism $G_k \cartprod K_2$, where the vertex in $G_k^1$ and $G_k^2$ corresponding to the vertex $v_j$ in $G_k$ is labeled $x_j$ and $y_j$, respectively, for $j \in [5k]$. Thus, the set $\cup_{j=1}^{5k} \{x_jy_j\}$ of edges between $V(G_k^1)$ and $V(G_k^2)$ in $G_k \cartprod K_2$ forms a perfect matching in $G_k \cartprod K_2$. For $i \in [k]$, let
\[
V_i = \bigcup_{j = 1}^5 \{x_{5(i-1)+j}, y_{5(i-1)+j} \}.
\]

When $k = 6$, the prism $G_k \cartprod K_2$ is illustrated in Figure~\ref{fig:prismG6}, where the vertices in $V_1$ are labelled. Let $S$ be an arbitrary total dominating set of $G_k \cartprod K_2$. For $i \in [k]$, let $S_i = S \cap V_i$. For $i \in [k]$, let
\[
X_i = \bigcup_{j = 2}^4 \{x_{5(i-1)+j}\} \quad \mbox{and} \quad
Y_i = \bigcup_{j = 2}^4 \{y_{5(i-1)+j}\}
\]

In order to totally dominate the vertices in the set $X_i$, we note that $|S_i| \ge 2$ for all $i \in [k]$.  Suppose that $|S_i| = 2$ for some $i \in [k]$. If both vertices in $S_i$ belong to the same copy of $G_k$, say to $G_k^2$, then at least one vertex in $X_i$ is not totally dominated by $S$. If the vertices in $S_i$ belong to different copies of $G_k$, then at least two vertices in $X_i \cup Y_i$ are not totally dominated by $S$. Both cases produce a contradiction, implying that $|S_i| \ge 3$. Hence,
\[
|S| = \sum_{i=1}^k |S_i| \ge 3k.
\]
Since $S$ is an arbitrary total dominating set of $G_k \cartprod K_2$, this implies that $\gt(G_k \cartprod K_2) \ge 3k$. To prove the converse, let
\[
X = \bigcup_{i=1}^{\lfloor k/2 \rfloor } \{ x_{10(i-1) + 1}, x_{10(i-1) + 2}, x_{10i} \}
\quad \mbox{and} \quad
Y = \bigcup_{i=1}^{\lfloor k/2 \rfloor} \{ y_{10(i-1) + 5}, y_{10(i-1) + 6},y_{10(i-1) + 7} \}.
\]
If $k$ is even, let
\[
D = (X \cup Y \cup \{x_{5k-1}\}) \setminus \{y_{5k-3}\}.
\]
If $k$ is odd, let
\[
D = X \cup Y \cup \{x_{5k-4},y_{5k-1},y_{5k}\}.
\]

For $k = 6$, the set $D$ is illustrated by the darkened vertices in Figure~\ref{fig:prismG6}. In both cases, $D$ is a total dominating set of $G_k \cartprod K_2$, and $|D \cap V_i| = 3$ for each $i \in [k]$, implying that
\[
\gt(G_k \cartprod K_2) \le |D| = \sum_{i=1}^k |D \cap V_i| = 3k.
\]
Consequently, $\gt(G_k \cartprod K_2) = 3k$.~\smallqed 

\medskip
By Claim~A and Claim~B, for $k \ge 2$, $\gamma(G_k) = 2k$ and  $\gt(G_k \cartprod K_2) = 3k$. This completes the proof of Theorem~\ref{t:nonbipartite}.
\qed

\begin{figure}[!ht]
\centering
\begin{tikzpicture}[scale=0.8,style=thick]
\def\vr{3pt}
\def\len{1}

\coordinate(x1) at (0,0); \coordinate(x2) at (0.5,0.7); \coordinate(x3) at (1.0,0); \coordinate(x4) at (1.5,0.7); \coordinate(x5) at (2.0,0);
\coordinate(y1) at (0,2); \coordinate(y2) at (0.5,2.7); \coordinate(y3) at (1.0,2); \coordinate(y4) at (1.5,2.7); \coordinate(y5) at (2.0,2);
\coordinate(x6) at (3,0); \coordinate(x7) at (3.5,0.7); \coordinate(x8) at (4.0,0); \coordinate(x9) at (4.5,0.7); \coordinate(x10) at (5.0,0);
\coordinate(y6) at (3,2); \coordinate(y7) at (3.5,2.7); \coordinate(y8) at (4.0,2); \coordinate(y9) at (4.5,2.7); \coordinate(y10) at (5.0,2);
\coordinate(x11) at (6,0); \coordinate(x12) at (6.5,0.7); \coordinate(x13) at (7.0,0); \coordinate(x14) at (7.5,0.7); \coordinate(x15) at (8.0,0);
\coordinate(y11) at (6,2); \coordinate(y12) at (6.5,2.7); \coordinate(y13) at (7.0,2); \coordinate(y14) at (7.5,2.7); \coordinate(y15) at (8.0,2);
\coordinate(x16) at (9,0); \coordinate(x17) at (9.5,0.7); \coordinate(x18) at (10.0,0); \coordinate(x19) at (10.5,0.7); \coordinate(x20) at (11.0,0);
\coordinate(y16) at (9,2); \coordinate(y17) at (9.5,2.7); \coordinate(y18) at (10.0,2); \coordinate(y19) at (10.5,2.7); \coordinate(y20) at (11.0,2);
\coordinate(x21) at (12,0); \coordinate(x22) at (12.5,0.7); \coordinate(x23) at (13.0,0); \coordinate(x24) at (13.5,0.7); \coordinate(x25) at (14.0,0);
\coordinate(y21) at (12,2); \coordinate(y22) at (12.5,2.7); \coordinate(y23) at (13.0,2); \coordinate(y24) at (13.5,2.7); \coordinate(y25) at (14.0,2);
\coordinate(x26) at (15,0); \coordinate(x27) at (15.5,0.7); \coordinate(x28) at (16.0,0); \coordinate(x29) at (16.5,0.7); \coordinate(x30) at (17.0,0);
\coordinate(y26) at (15,2); \coordinate(y27) at (15.5,2.7); \coordinate(y28) at (16.0,2); \coordinate(y29) at (16.5,2.7); \coordinate(y30) at (17.0,2);

\draw (x1) -- (x2) -- (x4) -- (x5) -- (x3) -- (x1);
\draw (y1) -- (y2) -- (y4) -- (y5) -- (y3) -- (y1);
\draw (x6) -- (x7) -- (x9) -- (x10) -- (x8) -- (x6);
\draw (y6) -- (y7) -- (y9) -- (y10) -- (y8) -- (y6);
\draw (x11) -- (x12) -- (x14) -- (x15) -- (x13) -- (x11);
\draw (y11) -- (y12) -- (y14) -- (y15) -- (y13) -- (y11);
\draw (x16) -- (x17) -- (x19) -- (x20) -- (x18) -- (x16);
\draw (y16) -- (y17) -- (y19) -- (y20) -- (y18) -- (y16);
\draw (x21) -- (x22) -- (x24) -- (x25) -- (x23) -- (x21);
\draw (y21) -- (y22) -- (y24) -- (y25) -- (y23) -- (y21);
\draw (x26) -- (x27) -- (x29) -- (x30) -- (x28) -- (x26);
\draw (y26) -- (y27) -- (y29) -- (y30) -- (y28) -- (y26);
\foreach \i in {1,...,30}
  { \draw (x\i) -- (y\i); }
\draw (x5) -- (x6); \draw (y5) -- (y6);
\draw (x10) -- (x11); \draw (y10) -- (y11);
\draw (x15) -- (x16); \draw (y15) -- (y16);
\draw (x20) -- (x21); \draw (y20) -- (y21);
\draw (x25) -- (x26); \draw (y25) -- (y26);

\foreach \i in {1,...,30}
{ \draw(x\i)[fill=white] circle(\vr); }
\foreach \i in {1,...,30}
{ \draw(y\i)[fill=white] circle(\vr); }
\draw(x1)[fill=black] circle(\vr);
\draw(x2)[fill=black] circle(\vr);
\draw(y5)[fill=black] circle(\vr);
\draw(y6)[fill=black] circle(\vr);
\draw(y7)[fill=black] circle(\vr);
\draw(x10)[fill=black] circle(\vr);
\draw(x11)[fill=black] circle(\vr);
\draw(x12)[fill=black] circle(\vr);
\draw(y15)[fill=black] circle(\vr);
\draw(y16)[fill=black] circle(\vr);
\draw(y17)[fill=black] circle(\vr);
\draw(x20)[fill=black] circle(\vr);
\draw(x21)[fill=black] circle(\vr);
\draw(x22)[fill=black] circle(\vr);
\draw(y25)[fill=black] circle(\vr);
\draw(y26)[fill=black] circle(\vr);
\draw(x29)[fill=black] circle(\vr);
\draw(x30)[fill=black] circle(\vr);

\draw(x1)node[below] {{\small $x_1$}};
\draw (0.6,0.35) node {{\small $x_2$}};
\draw(x3)node[below] {{\small $x_3$}};
\draw (1.4,0.35) node {{\small $x_4$}};
\draw(x5)node[below] {{\small $x_5$}};
\draw(-0.1,2.3) node {{\small $y_1$}};
\draw(y2)node[above] {{\small $y_2$}};
\draw(y3)node[above] {{\small $y_3$}};
\draw(y4)node[above] {{\small $y_4$}};
\draw(2.1,2.3) node {{\small $y_5$}};

\end{tikzpicture}

\caption{The prism $G_6 \cartprod K_2$}
\label{fig:prismG6}
\end{figure}

\section{Concluding Remarks}

Let us say that a graph $G$ is {\em $\gamma_t$-prism perfect} if $\gamma_t(G\cartprod K_2) = 2\gamma(G)$. We have seen that all bipartite graphs are $\gamma_t$-prism perfect. It would certainly be interesting to characterize $\gamma_t$-prism perfect graphs in general, but this appears to be a challenging problem. Instead, one could try to characterize $\gamma_t$-prism perfect graphs within some interesting families of graphs, say triangle-free graphs. 

A computation shows that among the $11.117$ connected graphs of order $8$, precisely $297$ graphs are not $\gamma_t$-prism perfect. Similarly, there are $79.638$ graphs that are not $\gamma_t$-prism perfect among the $11.716.571$ connected graphs of order $9$. These computations led us to conjecture the following conjecture. 

\begin{conj}
Almost all graphs are $\gamma_t$-prism perfect.
\end{conj}

With respect to the conjecture we refer to~\cite{glebov-2015} for the investigation of the behavior of the domination number in random graphs. 

Motivated by the construction presented in the proof of Theorem~\ref{t:nonbipartite} we wonder whether the following lower bound on the total domination number of prisms holds true. If so, then the construction implies that the bound is sharp. 

\begin{problem}
Is it true that for any graph $G$, $\gamma_t(G\cartprod K_2) \ge \frac{3}{2}\gamma(G)$?
\end{problem}

One may be tempted to try to extend the presented results to additional Cartesian product graphs. Clearly, $\gamma(P_3) = 2$ and an easy computation gives $\gamma_t(P_3 \square K_3) = \gamma_t(P_3 \square P_3) = 4$. Similarly $\gamma_t(P_3 \square K_4) = 4$ and $\gamma_t( P_3 \square P_4) = 6$, indicating that our result cannot be extended by a matter of parity. Moreover for all listed Cartesian products we were able to find pairs of bipartite graphs with the same domination number so that the total domination number of the respective Cartesian product differs. These examples give a strong evidence that the identity of Theorem~\ref{t:prism} cannot be generalized in ``obvious" directions. 

\section*{Acknowledgements}

Research of M.A.H.\ is supported in part by the South African National Research Foundation and the University of Johannesburg. J.A. and S.K. are supported by the Ministry of Science of Slovenia under the grants P1-0297 and N1-0043.

\end{document}